# SELF-INTERSECTION LOCAL TIME: CRITICAL EXPONENT, LARGE DEVIATIONS, AND LAWS OF THE ITERATED LOGARITHM[1]

By Richard F. Bass and Xia Chen

*University of Connecticut and University of Tennessee*

If $\beta_t$ is renormalized self-intersection local time for planar Brownian motion, we characterize when $\mathbb{E}e^{\gamma\beta_1}$ is finite or infinite in terms of the best constant of a Gagliardo–Nirenberg inequality. We prove large deviation estimates for $\beta_1$ and $-\beta_1$. We establish lim sup and lim inf laws of the iterated logarithm for $\beta_t$ as $t \to \infty$.

**1. Introduction.** Let $\beta_t$ be the renormalized self-intersection local time of a planar Brownian motion $X_t$. Formally,

$$\beta_t = \int_0^t \int_0^s \delta_0(X_s - X_u) \, du \, ds - \mathbb{E}\int_0^t \int_0^s \delta_0(X_s - X_u) \, du \, ds,$$

where $\delta_0$ is the delta function, and more precisely,

$$(1.1) \quad \beta_t = \lim_{\varepsilon \to 0}\left[\int_0^t \int_0^s \varphi_\varepsilon(X_s - X_u) \, du \, ds - \mathbb{E}\int_0^t \int_0^s \varphi_\varepsilon(X_s - X_u) \, du \, ds\right],$$

where $\varphi_\varepsilon$ is a suitable approximation to the identity. We have three main results in this paper:

1. Le Gall [16] showed that there is a critical exponent $\gamma_\beta$ such that

$$(1.2) \qquad \mathbb{E}e^{\gamma\beta_1} \begin{cases} < \infty, & \text{if } \gamma < \gamma_\beta, \\ = \infty, & \text{if } \gamma > \gamma_\beta. \end{cases}$$

   We characterize $\gamma_\beta$ in terms of the best constant of one of the Gagliardo–Nirenberg inequalities.
2. We prove large deviation estimates for $\beta_1$ and $-\beta_1$.

Received February 2003; revised July 2003.
[1]Supported in part by NSF Grants DMS-99-88496 and DMS-01-02238.

3. We prove laws of the iterated logarithm for the lim sup and lim inf behavior of $\beta_t$.

Self-intersection local time has been an object of much study in recent years. We cite [3, 5, 12, 13, 15, 16, 18, 19, 20, 22, 23] as an incomplete list of publications on this subject. In addition to probability theory, self-intersection local time has applications to some branches of mathematical physics, for example, constructive quantum field theories and polymer measures.

The quantity $\int_0^t \int_0^s \varphi_\varepsilon(X_s - X_u)\, du\, ds$ converges almost surely to infinity as $\varepsilon \to 0$ and to get convergence, the expectation of this integral must be subtracted. Therefore, exponential integrability of $\beta_1$ is a subtle issue. In 1994 Le Gall [16] proved there is a critical value $\gamma_\beta$ such that (1.2) holds. This fact has proved to be of considerable interest to the study of constructive quantum field theories. See also Theorem 2.23 of [5] for a discussion in the context of random walks with continuous time but discrete values. Our first main result characterizes $\gamma_\beta$.

THEOREM 1.1. *We have $\gamma_\beta = A^{-4}$, where $A > 0$ is the best constant in the inequality*

$$(1.3) \qquad \|f\|_4 \leq C\sqrt{\|\nabla f\|_2}\sqrt{\|f\|_2}, \qquad f:\mathbb{R}^2 \to \mathbb{R}.$$

Inequality (1.3) is one of a class of inequalities known as Gagliardo–Nirenberg inequalities. The proof of (1.3) is quite simple. Begin with the well-known Sobolev inequality in $\mathbb{R}^2$:

$$\|g\|_2 \leq c_1 \|\nabla g\|_1.$$

Replace $g$ by $f^2$, write $\nabla f^2$ as $2f \nabla f$ and apply the Cauchy–Schwarz inequality to the right-hand side. The best constant in (1.3) appears to be a difficult problem, however, and is currently open. The best constant for Nash's inequality, which is another special case of the Gagliardo–Nirenberg inequalities, was found by Carlen and Loss [6]. Two recent articles [9, 10] found the best constants for a class of Gagliardo–Nirenberg inequalities. Numerical values for the best constant in (1.3) were investigated as long ago as 1983 by Weinstein [21], who solved an eigenvalue problem by numerical methods and found that $A$ is approximately $(\pi \times 1.86225\ldots)^{-1/4}$. By Theorem 1.1,

$$\gamma_\beta \approx \pi \times 1.86225 \cdots \approx 5.85043.$$

This is very close to a conjecture made by B. Duplantier (private communication).

We could ask an analogous question about the intersection local time of two independent planar Brownian motions. There is a critical exponent $\gamma_\alpha$.

The critical value in this case was determined in [7] and was found to be the same constant $A^{-4}$ with $A$ as above. As a matter of fact, the result given in [7] is an important ingredient in the proof of Theorem 1.1 (and Theorem 1.2 as well).

As part of our proof of Theorem 1.1, we obtain large deviation estimates for $\beta_1$.

THEOREM 1.2. *We have*
$$\lim_{t \to \infty} \frac{1}{t} \log \mathbb{P}(\beta_1 \geq t) = -A^{-4},$$
*where $A$ is as in the statement of Theorem* 1.1.

We easily see that Theorem 1.1 is a direct consequence of Theorem 1.2. Interestingly, the lower tail of $\beta_1$ is not exponential, but instead is double exponential.

THEOREM 1.3. *There exists $0 < L < \infty$ such that*
$$\lim_{t \to \infty} t^{-2\pi} \log \mathbb{P}\{-\beta_1 \geq \log t\} = -L.$$

We also investigate laws of the iterated logarithm for $\beta_t$.

THEOREM 1.4. *We have*
$$\limsup_{t \to \infty} \frac{\beta_t}{t \log \log t} = \frac{1}{\gamma_\beta} \qquad a.s.$$

The lim inf behavior is described by the following theorem.

THEOREM 1.5. *We have*
$$\liminf_{t \to \infty} \frac{\beta_t}{t \log \log \log t} = -\frac{1}{2\pi} \qquad a.s.$$

Note the triple log in the rate of growth of the lim inf. This is suggested by the double exponential tail of $-\beta_t$. Compare this also with the result in [4] on the law of the iterated logarithm for the range of a random walk on $\mathbb{Z}^2$; the rate of growth there also has a triple log term. For a random walk the number of self-intersections is related to the range of the random walk up to time $n$, and Theorem 1.5 may provide some further insight into the result in [4]. Theorem 1.5 suggests that the right constant in [4] should be related to $1/2\pi$; we hope to return to these matters in future research.

Section 2 contains some basic facts about intersection local time. Theorems 1.1–1.3 are proved in Section 3. Theorem 1.4 is proved in Section 4 and Theorem 1.5 is proved in Section 5.

In all of the proofs, a key step is the representation of $\beta$ as the normalized sum of intersection local times of various pieces of the Brownian path plus sums of self-intersection local times; see Proposition 2.2. What makes the two-dimensional case much more difficult than the three-dimensional case is that in two dimensions these intersection local times of distinct pieces of the Brownian path are the dominant term.

**2. Preliminaries.** Let us begin with some notation. Let $X_t$ be a planar Brownian motion, let $\mathcal{F}_t$ be the completion of $\sigma\{X_s; s \leq t\}$ and let $\mathbb{P}^x$ denote the law of $X$ when $X$ is started at $x$. We use $\mathbb{P}$ for $\mathbb{P}^0$. The shift operators are denoted by $\theta_t$ as usual. If $I$ is an interval, we write $X(I)$ for the random set $\{X_s; s \in I\}$. The ball of radius $r$ about $x$ is denoted $B(x,r)$ and the letter $c$ with subscripts is used for positive finite constants whose exact value is unimportant.

If $X$ and $Y$ are two independent planar Brownian motions, the intersection local time can be defined formally by

$$\alpha(s,t) = \int_0^s \int_0^t \delta_0(X_r - Y_u)\,du\,dr,$$

where $\delta_0$ is the delta function. To make this rigorous, let $\varphi$ be a smooth nonnegative function in the Schwartz class which integrates to 1, let $\varphi_\varepsilon(x) = \varepsilon^{-2}\varphi(x/\varepsilon)$ (so that $\varphi_\varepsilon$ is an approximation to the identity) and define

$$(2.1) \qquad \alpha(s,t) = \lim_{\varepsilon \to 0} \int_0^s \int_0^t \varphi_\varepsilon(X_r - Y_u)\,du\,dr.$$

On the other hand, self-intersection local time cannot be defined so simply because the limit

$$\lim_{\varepsilon \to 0} \int_0^s \int_0^t \varphi_\varepsilon(X_r - X_u)$$

does not exist. A procedure called renormalization is needed. The renormalized self-intersection local time of $X$ is formally defined as

$$\beta_t = \int_0^t \int_0^s \delta_0(X_s - X_u)\,du\,ds - \mathbb{E}\int_0^t \int_0^s \delta_0(X_s - X_u)\,du\,ds.$$

To give a rigorous definition, let

$$(2.2) \quad \beta_t = \lim_{\varepsilon \to 0}\bigg[\int_0^t \int_0^s \varphi_\varepsilon(X_s - X_u)\,du\,ds - \mathbb{E}\int_0^t \int_0^s \varphi_\varepsilon(X_s - X_u)\,du\,ds\bigg].$$

That the limit exists a.s. and is continuous in $t$ is proved, for instance, in [13, 15] and [23]. Sometimes slightly different normalizations are used; they differ from ours by at most a constant times $t$. So there is no difference in the critical exponent or laws of the iterated logarithm, no matter which normalization is used.

If $I$ is an interval, we use $B(I)$ for the renormalized self-intersection for the piece of the path $X(I)$. That is, if $I = [s,t]$, then

$$B(I) = \beta_{t-s} \circ \theta_s. \tag{2.3}$$

If $I$ and $J$ are two intervals whose interiors are disjoint, let $A(I;J)$ denote the intersection local time for the two processes $X(I)$ and $X(J)$. To define this more precisely,

$$A(I;J) = \lim_{\varepsilon \to 0} \int_I \int_J \varphi_\varepsilon(X_s - X_t) \, ds \, dt. \tag{2.4}$$

REMARK 2.1. It is immediate by Brownian scaling that $\alpha(t,t)$ is equal in law to $t\alpha(1,1)$ and $\beta_t$ is equal in law to $t\beta_1$. Suppose $I = [a,b]$ and $J = [b,c]$. Then $A(I;J)$ measures the intersections of the two independent Brownian motions $X_b - X_{b-s}$ and $X_{b+t} - X_b$, and so $A(I;J)$ is equal in law to $\alpha(b-a, c-b)$ with starting point $(0,0)$.

PROPOSITION 2.2. *If $I$ is an interval that is the union of subintervals $I_j$, $1 \leq j \leq n$, such that the interiors of the $I_j$ are pairwise disjoint, then*

$$B(I) = \sum_{j=1}^n B(I_j) + \sum_{i<j} A(I_i; I_j) - \mathbb{E} \sum_{i<j} A(I_i; I_j). \tag{2.5}$$

PROOF. We have
$$\iint_{s,t \in I, s<t} \varphi_\varepsilon(X_t - X_s) \, ds \, dt - \mathbb{E} \iint_{s,t \in I, s<t} \varphi_\varepsilon(X_t - X_s) \, ds \, dt$$
$$= \sum_{j=1}^n \left[ \iint_{s,t \in I_j, s<t} \varphi_\varepsilon(X_t - X_s) \, ds \, dt - \mathbb{E} \iint_{s,t \in I_j, s<t} \varphi_\varepsilon(X_t - X_s) \, ds \, dt \right]$$
$$+ \sum_{i<j} \iint_{s \in I_i, t \in I_j} \varphi_\varepsilon(X_t - X_s) \, ds \, dt$$
$$- \sum_{i<j} \mathbb{E} \iint_{s \in I_i, t \in I_j} \varphi_\varepsilon(X_t - X_s) \, ds \, dt.$$
We now let $\varepsilon \to 0$. □

Although $\mathbb{E}\alpha(t,t)$ is a constant times $t$, we need a bit more precision.

PROPOSITION 2.3. *Let $\mathbb{P}^{(x_0,y_0)}$ be the joint law of $(X_t, Y_t)$ when $X_t$ is started at $x_0$ and $Y_t$ is started at $y_0$. Then*

$$\mathbb{E}^{(x_0,y_0)} \alpha(s,t) \leq \frac{1}{2\pi}[(s+t)\log(s+t) - s\log s - t\log t]. \tag{2.6}$$

*If $x_0 = y_0$, then we have equality in* (2.6).

PROOF. We have that $X_r$ is a two-dimensional normal random vector with mean $x_0$ and covariance matrix that is $r$ times the identity and that $Y_u$ is a two-dimensional normal random vector with mean $y_0$ and covariance matrix that is $u$ times the identity; moreover, the two random vectors are independent. Therefore, $X_r - Y_u$ is a two-dimensional normal random vector with mean $x_0 - y_0$ and covariance matrix that is $r + u$ times the identity. Hence

$$\mathbb{E}^{(x_0,y_0)} \int_0^s \int_0^t \varphi_\varepsilon(X_r - Y_u) \, dr \, du$$
$$= \int_0^s \int_0^t \int_{\mathbb{R}^2} \varphi_\varepsilon(z) \frac{1}{2\pi(r+u)} \exp\left(\frac{-|z - x_0 + y_0|^2}{2(r+u)}\right) dz \, dr \, du.$$

Letting $\varepsilon \to 0$ and using (2.1),

$$\mathbb{E}\alpha(s,t) = \int_0^s \int_0^t \frac{1}{2\pi(r+u)} \exp\left(\frac{-|x_0 - y_0|^2}{2(r+u)}\right) dr \, du.$$

The right-hand side is less than or equal to

$$\int_0^s \int_0^t \frac{1}{2\pi(r+u)} \, dr \, du$$

with equality when $x_0 = y_0$. Some routine calculus completes the proof. □

Le Gall [16] showed that there exists a value $\gamma_\beta$ such that

(2.7) $$\mathbb{E}e^{\gamma\beta_1} \begin{cases} < \infty, & \text{if } \gamma < \gamma_\beta, \\ = \infty, & \text{if } \gamma > \gamma_\beta. \end{cases}$$

In the same article, Le Gall proved that there exists a value $\gamma_\alpha$ such that

(2.8) $$\mathbb{E}e^{\gamma\alpha(1,1)} \begin{cases} < \infty, & \text{if } \gamma < \gamma_\alpha, \\ = \infty, & \text{if } \gamma > \gamma_\alpha. \end{cases}$$

He also gave a proof ([16], page 178) of a result by Varadhan [20] that

(2.9) $$\mathbb{E}e^{-\gamma\beta_1} < \infty$$

for all $\gamma > 0$.

**3. Large deviation estimates.** In [7], the large deviations for intersection local time of $p$ independent $d$-dimensional Brownian motions under the condition $p(d-2) < d$ were studied. Taking $d = p = 2$ in this result,

(3.1) $$\lim_{t\to\infty} \frac{1}{t} \log \mathbb{P}\{\alpha(1,1) \geq t\} = -A^{-4},$$

where $A > 0$ is the best constant in the Gagliardo–Nirenberg inequality

$$\|f\|_4 \leq C\sqrt{\|\nabla f\|_2}\sqrt{\|f\|_2}.$$

Let
$$M = \sup_{f \in \mathbb{F}_2} \left\{ \left( \int_{\mathbb{R}^2} |f(x)|^4 \, dx \right)^{1/2} - \tfrac{1}{2} \int_{\mathbb{R}^2} |\nabla f(x)|^2 \, dx \right\},$$
where $\mathbb{F}_2$ is the set of absolutely continuous functions on $\mathbb{R}^2$ satisfying
$$\int_{\mathbb{R}^2} |f(x)|^2 \, dx = 1 \quad \text{and} \quad \int_{\mathbb{R}^2} |\nabla f(x)|^2 \, dx < \infty.$$
As a special case of Lemma 8.2 in [7],
$$(3.2) \qquad M = \tfrac{1}{2} A^4.$$

In the following result, we claim that $\beta_1$ satisfies the same large deviation principle that $\alpha(1,1)$ does.

THEOREM 3.1.
$$(3.3) \qquad \lim_{t \to \infty} \frac{1}{t} \log \mathbb{P}\{\beta_1 \geq t\} = -A^{-4}.$$

*In particular,*
$$\mathbb{E} e^{\gamma \beta_1} \begin{cases} < \infty, & \text{if } \gamma < A^{-4}, \\ = \infty, & \text{if } \gamma > A^{-4}. \end{cases}$$

Note that this theorem implies $\gamma_\beta = A^{-4}$ and is a reformulation of Theorems 1.1 and 1.2.

PROOF. To establish the upper bound, we consider the decomposition
$$\beta_1 = \beta_{1/2} + B([1/2, 1]) + A([0, 1/2]; [1/2, 1]) - \mathbb{E} A([0, 1/2]; [1/2, 1]).$$
Recall that $B([1/2, 1])$ and $\beta_{1/2}$ are equal in law to $\tfrac{1}{2} \beta_1$, and $A([0, 1/2]; [1/2, 1])$ is equal in law to $\tfrac{1}{2} \alpha(1, 1)$. Moreover, $B([1/2, 1])$ is independent of $\beta_{1/2}$. Given $\varepsilon > 0$,
$$\mathbb{P}\{\beta_1 \geq t\} \leq \mathbb{P}\{\alpha(1,1) - \mathbb{E}\alpha(1,1) \geq (1-\varepsilon)t\} + \mathbb{P}\{\beta_1 + \beta_1' \geq (1+\varepsilon)t\},$$
where $\beta_1'$ is an independent copy of $\beta_1$. In view of (3.1),
$$\limsup_{t \to \infty} \frac{1}{t} \log \mathbb{P}\{\beta_1 \geq t\}$$
$$\leq \max\left\{ -(1-\varepsilon) A^{-4}, \limsup_{t \to \infty} \frac{1}{t} \log \mathbb{P}\{\beta_1 + \beta_1' \geq (1+\varepsilon)t\} \right\}.$$
We now need the simple fact that (1.2) is equivalent to
$$\limsup_{t \to \infty} \frac{1}{t} \log \mathbb{P}\{\beta_1 \geq t\} = -\gamma_\beta.$$

Also notice that

$$\mathbb{E}\exp\{\gamma(\beta_1 + \beta_1')\} = (\mathbb{E}\exp\{\gamma\beta_1\})^2 \begin{cases} < \infty, & \gamma < \gamma_\beta, \\ = \infty, & \gamma > \gamma_\beta. \end{cases}$$

So

$$\limsup_{t\to\infty} \frac{1}{t}\log \mathbb{P}\{\beta_1 + \beta_1' \geq (1+\varepsilon)t\} = -(1+\varepsilon)\gamma_\beta$$

and therefore

$$\limsup_{t\to\infty} \frac{1}{t}\log \mathbb{P}\{\beta_1 \geq t\} \leq -(1-\varepsilon)A^{-4}.$$

Letting $\varepsilon \to 0^+$, we obtain

(3.4) $$\limsup_{t\to\infty} \frac{1}{t}\log \mathbb{P}\{\beta_1 \geq t\} \leq -A^{-4}.$$

By scaling we have the upper bound of (3.3).

By scaling, Theorem 3.1 is equivalent to

(3.5) $$\lim_{n\to\infty} \frac{1}{n}\log \mathbb{P}\{\beta_n \geq \theta n^2\} = -\theta A^{-4}, \qquad \theta > 0.$$

Let

$$C_n = \sum_{k=1}^{n-1} A([0,k];[k,k+1]), \qquad n = 1, 2, \dots.$$

Then by Proposition 2.2,

$$\beta_n = C_n - \mathbb{E}C_n + \sum_{k=1}^{n} \beta([k-1,k]).$$

Notice that $\{\beta([k-1,k])\}$ is an i.i.d. sequence with the same distribution as $\beta_1$. Since the moment generating function of $\beta_1$ exists in a neighborhood of the origin, Cramér's theorem implies that for any $\delta > 0$,

(3.6) $$\lim_{n\to\infty} \frac{1}{n}\log \mathbb{P}\left\{\sum_{k=1}^{n} \beta([k-1,k]) \geq \delta n^2\right\} = -\infty.$$

Also, using Proposition 2.3, a calculation implies

(3.7) $$\mathbb{E}C_n = \frac{1}{2\pi} n\log n.$$

By Theorem 4.2.13 in [11], (3.5) is then equivalent to

(3.8) $$\lim_{n\to\infty} \frac{1}{n}\log \mathbb{P}\{C_n \geq \theta n^2\} = -\theta A^{-4}, \qquad \theta > 0.$$

We now claim that Theorem 3.1 holds provided

(3.9) $$\liminf_{n\to\infty} \frac{1}{n} \log \mathbb{E}\exp\{\lambda C_n^{1/2}\} \geq \frac{\lambda^2 A^2}{4}, \qquad \lambda > 0.$$

Indeed, from the upper bound (3.4), we can improve (3.9) into equality. In the case $\lambda < 0$, we use Jensen's inequality:

$$\mathbb{E}\exp\{\lambda C_n^{1/2}\} \geq \exp\{\lambda \mathbb{E} C_n^{1/2}\} \geq \exp\{\lambda (\mathbb{E} C_n)^{1/2}\} = \exp\{-O(\sqrt{n\log n})\},$$

where the last step follows from (3.7). Therefore, we have

$$\lim_{n\to\infty} \frac{1}{n} \log \mathbb{E}\exp\{\lambda C_n^{1/2}\} = \psi(\lambda)$$

for any real number $\lambda$, where

$$\psi(\lambda) = \begin{cases} \dfrac{\lambda^2 A^4}{4}, & \lambda \geq 0, \\ 0, & \lambda < 0. \end{cases}$$

By the Gärtner–Ellis theorem (Theorem 2.3.6 in [11]),

$$\lim_{n\to\infty} \frac{1}{n} \log \mathbb{P}\{C_n^{1/2} \geq \theta n\}$$
$$= -\sup_{\lambda\in\mathbb{R}}\{\lambda\theta - \psi(\lambda)\} = -\sup_{\lambda>0}\left\{\lambda\theta - \frac{\lambda^2 A^4}{4}\right\} = -\theta^2 A^{-4}, \qquad \theta > 0,$$

which is equivalent to (3.8).

We now prove (3.9). Some of the ideas come from [8]. We start with the fact (see, e.g., [17]) that for any measurable, bounded function $f$ on $\mathbb{R}^2$,

$$\lim_{n\to\infty} \frac{1}{n} \log \mathbb{E}\exp\left\{\int_0^n f(X_t)\,dt\right\} = \sup_{g\in\mathbb{F}_2}\left\{\int_{\mathbb{R}^2} f(x)g^2(x)\,dx - \frac{1}{2}\int_{\mathbb{R}^2} |\nabla g(x)|^2\,dx\right\}.$$

For any $\varepsilon > 0$, let $p_\varepsilon(x)$ be the density of $X_\varepsilon$ and write

$$L(t,x,\varepsilon) = \int_0^t p_\varepsilon(X_s - x)\,ds, \qquad x \in \mathbb{R}^2, t \geq 0.$$

It is easy to see from the semigroup property that

$$\left(\iint_{0\leq s\leq t\leq n} p_{2\varepsilon}(X_s - X_t)\,ds\,dt\right)^{1/2}$$
$$= \tfrac{1}{\sqrt{2}}\left(\int_{\mathbb{R}^2} L^2(n,x,\varepsilon)\,dx\right)^{1/2}$$
$$\geq \tfrac{1}{\sqrt{2}} \int_{\mathbb{R}^2} f(x) L(n,x,\varepsilon)\,dx = \tfrac{1}{\sqrt{2}} \int_0^n f_\varepsilon(X_t)\,dt$$

for any measurable $f$ on $\mathbb{R}^2$ with
$$\int_{\mathbb{R}^2} f^2(x)\,dx = 1,$$
where
$$f_\varepsilon(x) = \int_{\mathbb{R}^2} f(x-y) p_\varepsilon(y)\,dy.$$
Therefore,
$$\liminf_{n\to\infty} \frac{1}{n} \log \mathbb{E} \exp\left\{\lambda \left(\iint_{0\le s\le t\le n} p_{2\varepsilon}(X_s - X_t)\,ds\,dt\right)^{1/2}\right\}$$
$$\ge \sup_{g\in\mathbb{F}_2}\left\{\frac{\lambda}{\sqrt{2}} \int_{\mathbb{R}^2} f_\varepsilon(x) g^2(x)\,dx - \frac{1}{2}\int_{\mathbb{R}^2}|\nabla g(x)|^2\,dx\right\}$$
$$= \sup_{g\in\mathbb{F}_2}\left\{\frac{\lambda}{\sqrt{2}} \int_{\mathbb{R}^2} f(x)\left(\int_{\mathbb{R}^2} g^2(x-y) p_\varepsilon(y)\,dy\right)dx - \frac{1}{2}\int_{\mathbb{R}^2}|\nabla g(x)|^2\,dx\right\}.$$

Taking the supremum over $f$ with $\|f\|_2 = 1$ and using the fact that the dual of $L^2$ is $L^2$ gives

$$\liminf_{n\to\infty} \frac{1}{n} \log \mathbb{E} \exp\left\{\lambda \left(\iint_{0\le s\le t\le n} p_{2\varepsilon}(X_s - X_t)\,ds\,dt\right)^{1/2}\right\}$$
(3.10)
$$\ge \sup_{g\in\mathbb{F}_2}\left\{\frac{\lambda}{\sqrt{2}}\left[\int_{\mathbb{R}^2}\left(\int_{\mathbb{R}^2} g^2(x-y) p_\varepsilon(y)\,dy\right)^2 dx\right]^{1/2}\right.$$
$$\left. - \frac{1}{2}\int_{\mathbb{R}^2}|\nabla g(x)|^2\,dx\right\}$$

for any $\lambda > 0$.

On the other hand, write
$$\xi_k(\varepsilon) = \iint_{\{k-1\le s\le t\le k\}} p_{2\varepsilon}(X_s - X_t)\,ds\,dt, \qquad k = 1, 2, \ldots,$$
and
$$D_n = \bigcup_{k=1}^{n-1} [0, k) \times (k, k+1], \qquad n = 1, 2, \ldots.$$

Then $\{\xi_k(\varepsilon)\}_{k\ge 1}$ is an i.i.d. sequence and
$$\iint_{\{0\le s\le t\le n\}} p_{2\varepsilon}(X_s - X_t)\,ds\,dt = \iint_{D_n} p_{2\varepsilon}(X_s - X_t)\,ds\,dt + \sum_{k=1}^n \xi_k(\varepsilon).$$

Let $p, q > 1$ be such that $p^{-1} + q^{-1} = 1$. By the triangle inequality and Hölder's inequality,
$$\mathbb{E}\exp\left\{p^{-1}\lambda\left(\iint_{\{0\le s\le t\le n\}} p_{2\varepsilon}(X_s - X_t)\,ds\,dt\right)^{1/2}\right\}$$

$$\leq \left[\mathbb{E}\exp\left\{\lambda\left(\iint_{D_n} p_{2\varepsilon}(X_s - X_t)\,ds\,dt\right)^{1/2}\right\}\right]^{1/p}$$

$$\times \left[\mathbb{E}\exp\left\{qp^{-1}\lambda\left(\sum_{k=1}^{n}\xi_k(\varepsilon)\right)^{1/2}\right\}\right]^{1/q}.$$

It is easy to see from standard large deviation theory that

$$\lim_{n\to\infty} \frac{1}{n}\log \mathbb{E}\exp\left\{qp^{-1}\lambda\left(\sum_{k=1}^{n}\xi_k(\varepsilon)\right)^{1/2}\right\} = 0.$$

Therefore, by (3.10) we have

$$\liminf_{n\to\infty} \frac{1}{n}\log \mathbb{E}\exp\left\{\lambda\left(\iint_{D_n} p_{2\varepsilon}(X_s - X_t)\,ds\,dt\right)^{1/2}\right\}$$

$$\geq p\sup_{g\in\mathbb{F}_2}\left\{\frac{p^{-1}\lambda}{\sqrt{2}}\left[\int_{\mathbb{R}^2}\left(\int_{\mathbb{R}^2} g^2(x-y)p_\varepsilon(y)\,dy\right)^2 dx\right]^{1/2}\right.$$

$$\left. -\frac{1}{2}\int_{\mathbb{R}^2}|\nabla g(x)|^2\,dx\right\}.$$

Letting $p \to 1^+$ gives

$$\liminf_{n\to\infty} \frac{1}{n}\log \mathbb{E}\exp\left\{\lambda\left(\iint_{D_n} p_{2\varepsilon}(X_s - X_t)\,ds\,dt\right)^{1/2}\right\}$$

(3.11)
$$\geq \sup_{g\in\mathbb{F}_2}\left\{\frac{\lambda}{\sqrt{2}}\left[\int_{\mathbb{R}^2}\left(\int_{\mathbb{R}^2} g^2(x-y)p_\varepsilon(y)\,dy\right)^2 dx\right]^{1/2}\right.$$

$$\left. -\frac{1}{2}\int_{\mathbb{R}^2}|\nabla g(x)|^2\,dx\right\}.$$

For any $m \geq 0$, let $k \geq 0$ be the integer such that $2k \leq m \leq 2(k+1)$. By Lemma 3.4,

$$\mathbb{E}(C_n^{(m+2)/2}) \geq [\mathbb{E}C_n^{k+1}]^{(m+2)/(2(k+1))}$$

$$\geq \left[\mathbb{E}\left(\iint_{D_n} p_{2\varepsilon}(X_s - X_t)\,ds\,dt\right)^{k+1}\right]^{(m+2)/(2(k+1))}$$

$$\geq \left[\mathbb{E}\left(\iint_{D_n} p_{2\varepsilon}(X_s - X_t)\,ds\,dt\right)^{m/2}\right]^{(m+2)/m}.$$

As $n \to \infty$, it is clear that $C_n \to \infty$. Using Lemma 3.4, we can also see that there is a $N > 0$ and $\varepsilon_0 > 0$ such that

$$\mathbb{E}\left(\iint_{D_n} p_{2\varepsilon}(X_s - X_t)\,ds\,dt\right)^{m/2} \geq 1, \qquad m = 0, 1, \ldots,$$

if $n \geq N$ and $\varepsilon \leq \varepsilon_0$. Hence

$$\mathbb{E}(C_n^{(m+2)/2}) \geq \mathbb{E}\bigg(\iint_{D_n} p_{2\varepsilon}(X_s - X_t)\,ds\,dt\bigg)^{m/2}.$$

Using the Taylor series expansion for $e^{\lambda x}$, for each $0 < \delta < \lambda$,

$$\mathbb{E}(C_n \exp((\lambda - \delta)C_n^{1/2})) \geq \mathbb{E}\exp\bigg\{(\lambda - \delta)\bigg(\iint_{D_n} p_{2\varepsilon}(X_s - X_t)\,ds\,dt\bigg)^{1/2}\bigg\}.$$

By the fact that $e^{\lambda\sqrt{x}} \geq xe^{(\lambda-\delta)\sqrt{x}}$ for sufficiently large $x > 0$ and in view of (3.11) (with $\lambda$ replaced by $\lambda - \delta$), the above estimate implies

$$\liminf_{n\to\infty} \frac{1}{n} \log \mathbb{E}\exp\{\lambda C_n^{1/2}\}$$
$$\geq \sup_{g\in\mathbb{F}_2}\bigg\{\frac{\lambda-\delta}{\sqrt{2}}\bigg[\int_{\mathbb{R}^2}\bigg(\int_{\mathbb{R}^2} g^2(x-y)p_\varepsilon(y)\,dy\bigg)^2 dx\bigg]^{1/2} - \frac{1}{2}\int_{\mathbb{R}^2}|\nabla g(x)|^2\,dx\bigg\}.$$

Letting $\varepsilon \to 0^+$ on the right-hand side gives

(3.12)
$$\begin{aligned}\liminf_{n\to\infty} \frac{1}{n} &\log \mathbb{E}\exp\{\lambda C_n^{1/2}\} \\ &\geq \sup_{g\in\mathbb{F}_2}\bigg\{\frac{\lambda-\delta}{\sqrt{2}}\bigg(\int_{\mathbb{R}^2}|g(x)|^4\,dx\bigg)^{1/2} - \frac{1}{2}\int_{\mathbb{R}^2}|\nabla g(x)|^2\,dx\bigg\} \\ &= \frac{(\lambda-\delta)^2}{2}\sup_{f\in\mathbb{F}_2}\bigg\{\bigg(\int_{\mathbb{R}^2}|f(x)|^4\,dx\bigg)^{1/2} - \frac{1}{2}\int_{\mathbb{R}^2}|\nabla f(x)|^2\,dx\bigg\} \\ &= \frac{(\lambda-\delta)^2 A^4}{4}\end{aligned}$$

for any $0 < \delta < \lambda$, where the second step follows from the substitution

$$g(x) = \frac{\lambda-\delta}{\sqrt{2}} f\bigg(\frac{\lambda-\delta}{\sqrt{2}} x\bigg)$$

and the last step follows from (3.2). Finally, letting $\delta \to 0^+$ gives (3.9).  □

THEOREM 3.2. *There is a $0 < L \leq \infty$ such that*

$$\lim_{t\to\infty} t^{-2\pi} \log \mathbb{P}\{-\beta_1 \geq \log t\} = -L.$$

Theorem 3.2 proves part of Theorem 1.3.

PROOF OF THEOREM 3.2. For any positive integers $m$ and $n$, by Proposition 2.2,

$$\beta_{m+n} = \beta_n + B([n, n+m]) + A([0,n]; [n, n+m]) - \mathbb{E}A([0,n]; [n, n+m])$$
$$\geq \beta_n + B([n, n+m]) - \mathbb{E}A([0,n]; [n, n+m]),$$

and $\beta_n$ and $B([n, n+m])$ are independent. Hence

$$\mathbb{E}\exp\{-2\pi(m+n)\beta_1\}$$
$$= \mathbb{E}\exp\{-2\pi\beta_{m+n}\}$$
$$\leq \exp\{2\pi\mathbb{E}A([0,n];[n,n+m])\}\mathbb{E}\exp\{-2\pi B([n,n+m])\}\mathbb{E}\exp\{-2\pi\beta_n\}$$
$$= \exp\{2\pi\mathbb{E}A([0,n];[n,n+m])\}\mathbb{E}\exp\{-2\pi\beta_m\}\mathbb{E}\exp\{-2\pi\beta_n\}$$
$$= \exp\{2\pi\mathbb{E}A([0,n];[n,n+m])\}\mathbb{E}\exp\{-2\pi m\beta_1\}\mathbb{E}\exp\{-2\pi n\beta_1\}.$$

Using this and Proposition 2.3,

$$(m+n)^{-(m+n)}\mathbb{E}\exp\{-2\pi(m+n)\beta_1\}$$
$$\leq (m^{-m}\mathbb{E}\exp\{-2\pi m\beta_1\})(n^{-n}\mathbb{E}\exp\{-2\pi n\beta_1\}).$$

If we write

$$a(n) = \log(n^{-n}\mathbb{E}\exp\{-2\pi n\beta_1\}), \qquad n=1,2,\ldots,$$

then we have proved that for any positive integers $m$ and $n$,

$$a(n+m) \leq a(m) + a(n).$$

Consequently,

$$\lim_{n\to\infty}\frac{1}{n}a(n) = \inf_{m\geq 1}\left\{\frac{1}{m}a(m)\right\}.$$

By Stirling's formula, this is equivalent to

$$\lim_{n\to\infty}\frac{1}{n}\log((n!)^{-1}\mathbb{E}\exp\{-2\pi n\beta_1\}) = 1 + \inf_{m\geq 1}\left\{\frac{1}{m}a(m)\right\}.$$

By Lemma 2.3 of [14],

$$\limsup_{t\to\infty} t^{-1}\log\mathbb{P}\{\exp\{-2\pi\beta_1\} \geq t\} = -\exp\left\{-1 - \inf_{m\geq 1}\left\{\frac{1}{m}a(m)\right\}\right\} = -L,$$

where

$$L = \exp\left\{-1 - \inf_{m\geq 1}\left\{\frac{1}{m}a(m)\right\}\right\}. \qquad \square$$

REMARK 3.3. In fact, $L < \infty$. This is established in Corollary 5.7.

LEMMA 3.4. *For any positive numbers $\varepsilon$ and $\varepsilon'$ with $\varepsilon > \varepsilon'$, any $D \subset \{(s,t); s \leq t\}$ and integer $m \geq 1$,*

$$\mathbb{E}\left[\iint_D p_{\varepsilon'}(X_t - X_s)\,ds\,dt\right]^m \geq \mathbb{E}\left[\iint_D p_\varepsilon(X_t - X_s)\,ds\,dt\right]^m.$$

*Furthermore, if $D$ is a finite union of disjoint rectangles contained in $\{(s,t); s \leq t\}$,*

$$D = \bigcup_{k=1}^{n} (I_k \times J_k),$$

*then*

$$\mathbb{E}\left[\iint_D p_\varepsilon(X_t - X_s)\,ds\,dt\right]^m \leq \mathbb{E}\left[\sum_{k=1}^{n} A(I_k; J_k)\right]^m.$$

PROOF. By the Fourier transform,

$$p_\varepsilon(X_t - X_s) = \frac{1}{(2\pi)^2} \int_{\mathbb{R}^2} d\xi \exp\{-i\xi \cdot (X_t - X_s)\} \exp\left\{-\frac{\varepsilon}{2}|\xi|^2\right\}.$$

Hence

$$\mathbb{E}\left[\iint_D p_\varepsilon(X_t - X_s)\,ds\,dt\right]^m$$

$$= \frac{1}{(2\pi)^{2m}} \int_{D^m} ds_1\,dt_1 \cdots ds_m\,dt_m \int_{(\mathbb{R}^2)^m} d\xi_1 \cdots d\xi_m$$

$$\times \mathbb{E}\left[\prod_{k=1}^{m} \exp\{-i\xi_k \cdot (X_t - X_s)\} \exp\left\{-\frac{\varepsilon}{2}|\xi_k|^2\right\}\right]$$

$$= \frac{1}{(2\pi)^{2m}} \int_{D^m} ds_1\,dt_1 \cdots ds_m\,dt_m \int_{(\mathbb{R}^2)^m} d\xi_1 \cdots d\xi_m$$

$$\times \exp\left\{-\frac{1}{2}\mathrm{Var}\left[\sum_{k=1}^{m} \xi_k \cdot (X_t - X_s)\right]\right\} \exp\left\{-\frac{\varepsilon}{2}\sum_{k=1}^{m}|\xi_k|^2\right\},$$

which leads to the first half of the lemma.

As for the second half of the lemma, by Theorem 4 on page 191 in [15],

$$\mathbb{E}\left[\sum_{k=1}^{n} A(I_k; J_k)\right]^m < \infty, \qquad m = 0, 1, \ldots,$$

and

$$\iint_D p_\varepsilon(X_t - X_s)\,ds\,dt \to \sum_{k=1}^{n} A(I_k; J_k) \qquad (\varepsilon \to 0^+)$$

in $L^m$-norm for all integers $m \geq 1$. (In fact, Le Gall proved the above convergence with $p_\varepsilon$ replaced by the uniform density on the disk of radius $\varepsilon$. It can be seen from his argument that this remains true in our case.) Therefore, letting $\varepsilon' \to 0^+$ leads to the second half of the lemma. □

**4. The lim sup result.** In this section we establish Theorem 1.4.

LEMMA 4.1. *There exist constants $c_1, c_2$ such that for all $\lambda > 0$ and all $a \in (0,1)$,*
$$\mathbb{P}(\alpha(1,a) > \lambda) \leq c_1 \exp(-c_2 \lambda / \sqrt{a}).$$

PROOF. Let $m \geq 1$ be an integer. We first prove there exist constants $c_3, c_4$ such that
$$\mathbb{E}[\alpha(1,a)^m] \leq c_3 c_4^m a^{m/2} m!. \tag{4.1}$$

To establish this, write
$$\alpha(1,a) = \lim_{\varepsilon \to 0} \int_0^1 \int_0^a p(\varepsilon, 0, X_v - Y_u) \, dv \, du,$$

where $p(\varepsilon, x, y)$ is the transition density of planar Brownian motion. As mentioned in the last paragraph of the proof of Lemma 3.4, the convergence takes place in $L^p$ for every $p$. By the semigroup property,
$$p(\varepsilon, 0, X_v - Y_u) = p(\varepsilon, X_v, Y_u) = \int_{\mathbb{R}^2} p(\varepsilon/2, x, X_v) p(\varepsilon/2, x, Y_u) \, dx$$

and so
$$\left[ \int_0^1 \int_0^a p(\varepsilon, 0, X_v - Y_u) \, dv \, du \right]^m$$
$$= \int_{(\mathbb{R}^2)^m} dx_1 \cdots dx_m \left( \prod_{k=1}^m \int_0^1 p(\varepsilon/2, x_k, X_v) \, dv \right)$$
$$\times \left( \prod_{k=1}^m \int_0^a p(\varepsilon/2, x_k, Y_u) \, du \right).$$

Using the independence of $X$ and $Y$, the expectation is equal to
$$\int_{(\mathbb{R}^2)^m} dx_1 \cdots dx_m \, \mathbb{E}\left[ \prod_{k=1}^m \int_0^1 p(\varepsilon/2, x_k, X_v) \, dv \right] \mathbb{E}\left[ \prod_{k=1}^m \int_0^a p(\varepsilon/2, x_k, Y_u) \, du \right].$$

By the Cauchy–Schwarz inequality this is less than $J_1(\varepsilon)^{1/2} J_2(\varepsilon)^{1/2}$, where
$$J_1(\varepsilon) = \int_{(\mathbb{R}^2)^m} dx_1 \cdots dx_m \left( \mathbb{E}\left[ \prod_{k=1}^m \int_0^1 p(\varepsilon/2, x_k, X_v) \, dv \right] \right)^2$$

and
$$J_2(\varepsilon) = \int_{(\mathbb{R}^2)^m} dx_1 \cdots dx_m \left( \mathbb{E}\left[ \prod_{k=1}^m \int_0^a p(\varepsilon/2, x_k, Y_u) \, du \right] \right)^2.$$

By Brownian scaling,

(4.2) $$\lim_{\varepsilon \to 0} J_2(\varepsilon) = a^m \lim_{\varepsilon \to 0} J_1(\varepsilon).$$

To estimate $J_1(\varepsilon)$, we rewrite it as

$$\int_{(\mathbb{R}^2)^m} dx_1 \cdots dx_m \, \mathbb{E}\left[\prod_{k=1}^m \int_0^1 p(\varepsilon/2, x_k, X_v)\, dv\right] \mathbb{E}\left[\prod_{k=1}^m \int_0^1 p(\varepsilon/2, x_k, Y_u)\, du\right]$$

and so by the argument above in reverse order,

$$J_1(\varepsilon) = \mathbb{E}\left[\int_0^1 \int_0^1 p(\varepsilon, X_v, Y_u)\, dv\, du\right]^m.$$

Therefore $\lim_{\varepsilon \to 0} J_1(\varepsilon) = \mathbb{E}[\alpha(1,1)^m]$. Lemma 2 of [16] together with (4.2) and an application of Fatou's lemma completes the proof of (4.1).

We then obtain

$$\mathbb{E}\exp\left(\frac{\alpha(1,a)}{2c_4\sqrt{a}}\right) = \sum_{m=0}^\infty \left(\frac{1}{2c_4\sqrt{a}}\right)^m \frac{\mathbb{E}[\alpha(1,a)^m]}{m!} \leq c_7,$$

where $c_7$ does not depend on $a$. Finally,

$$\mathbb{P}(\alpha(1,a) > \lambda)$$
$$\leq \exp(-\lambda/(2c_4\sqrt{a}))\mathbb{E}\exp(\alpha(1,a)/(2c_4\sqrt{a})) \leq c_7 \exp(-\lambda/(2c_4\sqrt{a})),$$

which is what we wanted. $\square$

The key to the upper bound is to obtain an estimate of the following form.

PROPOSITION 4.2. *If $\gamma < \gamma_\beta$, there exists $c_1$ such that*

(4.3) $$\mathbb{P}\left(\sup_{t \leq 1} \beta_t > \lambda\right) \leq c_1 e^{-\gamma\lambda}, \qquad \lambda > 0.$$

PROOF. By Proposition 2.2,

$$\beta_t - \beta_s = B([s,t]) + A([0,s];[s,t]) - \mathbb{E}A([0,s];[s,t]).$$

Let $\gamma'$ be the midpoint of $(\gamma, \gamma_\beta)$ and let $\varepsilon > 0$ be chosen so that $\gamma'(1-\varepsilon)$ is the midpoint of $(\gamma, \gamma')$. Note

(4.4) $$\mathbb{P}(\beta_t - \beta_s > \lambda) \leq \mathbb{P}(B([s,t]) > \lambda/2) + \mathbb{P}(A([0,s];[s,t]) > \lambda/2).$$

Since $B([s,t])$ equals $\beta_{t-s}$ in law, which equals $(t-s)\beta_1$ in law, the first probability on the right is bounded by

(4.5) $$c_2 \exp\left\{-\frac{\gamma'\lambda}{2(t-s)}\right\}.$$

However, $A([0,s];[s,t])$ is equal in law to $\alpha(s,t-s)$, which is smaller than $\alpha(1,t-s)$. So by Lemma 4.1 there exists $c_3$ not depending on $s$ or $t$ such that

$$(4.6) \qquad \mathbb{P}\left(A([0,s];[s,t]) > \frac{\lambda}{2}\right) \leq \exp\left\{-\frac{c_3\lambda}{(t-s)^{1/2}}\right\}.$$

Fix $n = 2^N$. Since

$$\sup_{k \leq n} \mathbb{E}\exp(\gamma'\beta_{k/n}) = \sup_{k \leq n} \mathbb{E}\exp((\gamma'k/n)\beta_1) \leq c_4,$$

where $c_4$ does not depend on $n$, then

$$(4.7) \qquad \mathbb{P}\left(\sup_{k \leq n} \beta_{k/n} > (1-\varepsilon)\lambda\right) \leq n\sup_{k \leq n} \mathbb{P}(\beta_{k/n} > (1-\varepsilon)\lambda) \\ \leq ne^{-\gamma'(1-\varepsilon)\lambda} \leq ne^{-\gamma\lambda}.$$

Now we use metric entropy. If $t \in (0,1)$, let $t_j$ be the largest multiple of $2^{-j}$ that is less than or equal to $t$. Write

$$\beta_t = \beta_{t_N} + (\beta_{t_{N+1}} - \beta_{t_N}) + (\beta_{t_{N+2}} - \beta_{t_{N+1}}) + \cdots.$$

If $\beta_t > \lambda$ for some $t \leq 1$, either (a) for some $k \leq n$, we have $\beta_{k/n} > (1-\varepsilon)\lambda$ or (b) for some $j \geq N$ and some $s < t$ with $t - s = 2^{-j}$ and both $s,t$ integer multiples of $2^{-j}$, we have

$$(4.8) \qquad \beta_t - \beta_s > \varepsilon\lambda/(100j^2).$$

The probability of possibility (a) is bounded by (4.7). Using (4.5) and (4.6), the probability of possibility (b) is bounded by

$$(4.9) \qquad c_5 \sum_{j=N}^{\infty} 2^j[\exp(-\varepsilon\gamma'\lambda 2^j/(200j^2)) + \exp(-c_3\varepsilon\lambda 2^{j/2}/(200j^2))].$$

The $2^j$ in front of the brackets comes about because there are $2^j$ pairs $(s,t)$ to consider. It is not hard to see that the sum in (4.9) is bounded by

$$c_6[\exp(-\varepsilon\gamma'\lambda 2^N/(400N^2)) + \exp(-c_3\varepsilon\lambda 2^{N/2}/(400N^2))].$$

If we choose $N$ large enough so that $2^N\varepsilon/(400N^2) > 1$ and $c_3 2^{N/2}\varepsilon/(400 \times N^2) > \gamma'$, we then have that the probability of possibility (b) is bounded by

$$2c_7 e^{-\gamma'\lambda} \leq 2c_7 e^{-\gamma\lambda}.$$

If we combine this with (4.7), we have (4.3). $\square$

Using the Borel–Cantelli lemma it is now straightforward to get the following theorem:

THEOREM 4.3. *We have*
$$\limsup_{t\to\infty} \frac{\beta_t}{t\log\log t} \leq \frac{1}{\gamma_\beta} \qquad a.s.$$

PROOF. Let $M > 1/\gamma_\beta$. Choose $\varepsilon > 0$ small and $q > 1$ close to 1 so that $M(\gamma_\beta - 2\varepsilon)/q > 1$. Let $t_n = q^n$ and let $C_n = \{\sup_{s \leq t_n} \beta_s > Mt_{n-1}\log\log t_{n-1}\}$. By Proposition 4.2 and scaling, the probability of $C_n$ is bounded by
$$c_1 \exp(-(\gamma_\beta - \varepsilon)Mt_{n-1}\log\log t_{n-1}/t_n).$$
By our choices of $\varepsilon$ and $q$ this is summable, so by the Borel–Cantelli lemma the probability that $C_n$ happens infinitely often is zero. To complete the proof we point out that if $\beta_t > Mt\log\log t$ for some $t \in [t_{n-1}, t_n]$, then the event $C_n$ occurs. $\square$

To finish the proof of Theorem 4.3 we prove the next theorem:

THEOREM 4.4. *We have*
$$\limsup_{t\to\infty} \frac{\beta_t}{t\log\log t} \geq \frac{1}{\gamma_\beta} \qquad a.s.$$

Jointly, Theorems 4.3 and 4.4 are a reformulation of Theorem 1.4.

PROOF OF THEOREM 4.4. Let $\gamma > \gamma_\beta$ and let $\gamma'$ be the midpoint of $(\gamma_\beta, \gamma)$. Then by Theorem 3.1,

(4.10) $$\mathbb{P}(\beta_1 \geq a\log\log n) \geq c_2 e^{-\gamma' a\log\log n}, \qquad a > 0.$$

Let $\delta > 0$ be small enough so that $(1+\delta)\gamma'/\gamma < 1$ and set $t_n = \exp(n^{1+\delta})$. By (2.5),
$$\begin{aligned}\beta_{t_n} &= B([0, t_n]) \\ &= B([t_{n-1}, t_n]) + B([0, t_{n-1}]) \\ &\quad + A([0, t_{n-1}]; [t_{n-1}, t_n]) - \mathbb{E}A([0, t_{n-1}]; [t_{n-1}, t_n]) \\ &\geq B([t_{n-1}, t_n]) + B([0, t_{n-1}]) - \mathbb{E}A([0, t_{n-1}]; [t_{n-1}, t_n]).\end{aligned}$$

By scaling,
$$\begin{aligned}\mathbb{E}A([0, t_{n-1}]; [t_{n-1}, t_n]) &\leq \mathbb{E}\alpha(t_n, t_n) \\ &= t_n\mathbb{E}\alpha(1, 1) = o(t_n\log\log t_n), \qquad n\to\infty.\end{aligned}$$

Since $A \geq 0$, we need only to prove

(4.11) $$\limsup_{n\to\infty} \frac{B([t_{n-1}, t_n])}{t_n\log\log t_n} \geq \frac{1}{\gamma_\beta} \qquad a.s.$$

and

(4.12) $$\lim_{n\to\infty} \frac{|B([0,t_{n-1}])|}{t_n \log\log t_n} = 0 \quad \text{a.s.}$$

Using (4.10) and scaling, it is straightforward to obtain

$$\sum_{n=1}^{\infty} \mathbb{P}\left(B([t_{n-1},t_n]) > \frac{1}{\gamma} t_n \log\log t_n\right) = \infty.$$

Using the fact that different pieces of a Brownian path are independent and the Borel–Cantelli lemma,

$$\limsup_{n\to\infty} \frac{B([t_{n-1},t_n])}{t_n \log\log t_n} > \frac{1}{\gamma} \quad \text{a.s.}$$

Letting $\gamma \to \gamma_\beta^+$ gives (4.11).

Let $\varepsilon > 0$. By (2.9) there exists $c_4 > 0$ such that

(4.13) $$\mathbb{P}\{-\beta_1 \geq \varepsilon \log\log n\} \leq c_4 e^{-2\log\log n}.$$

So (4.12) follows from Theorem 3.1, (4.13), scaling and the Borel–Cantelli lemma. □

REMARK 4.5. Theorems 4.3 and 4.4 together imply Theorem 1.4.

**5. The lim inf result.** Let us write $D_t$ for $-\beta_t$. We know $\mathbb{E}\exp(\gamma D_1) < \infty$ for every $\gamma > 0$, but in fact we have the following proposition.

PROPOSITION 5.1. *We have $\mathbb{E}\exp(\gamma \sup_{s\leq 1} D_s) < \infty$ for every $\gamma > 0$.*

PROOF. Fix $\gamma > 0$. Choose $N$ a fixed integer so that $2^N/(1600N^3) > 2$. If $s < t \leq 1$, we know $\mathbb{E}A([0,s];[s,t]) \leq c_1(t-s)L(t-s) \leq c_2$, where $L(x) = 1 + |\log(1/x)|$. Suppose $\lambda > c_3$, where $c_3$ is chosen so that $c_3/(400j^2) > 2c_1 2^{-j} L(2^{-j})$ for each $j \geq 0$. Let $s_j = \inf\{k/2^j : s \leq k/2^j\}$. If $s \in [0,1]$, we can write

$$D_s = D_{s_N} + (D_{s_{N+1}} - D_{s_N}) + (D_{s_{N+2}} - D_{s_{N+1}}) + \cdots.$$

So if $D_s > \lambda$ for some $s \in [0,1]$, then either (a) for some $k \leq 2^N$, we have $D_{k/2^N} > \lambda/2$ or (b) for some $j > N$ and some $s < t$, both multiples of $2^{-j}$ with $t - s = 2^{-j}$,

$$D_t - D_s > \frac{\lambda}{200j^2}.$$

We have $\mathbb{P}(D_{k/2^N} > \lambda/2) = \mathbb{P}(D_1 > 2^N\lambda/(2k)) \leq \mathbb{P}(D_1 > \lambda/2) \leq c_4 e^{-2\gamma\lambda}$ since $\mathbb{E}\exp(4\gamma D_1) < \infty$. So the probability of possibility (a) is bounded by

$$c_4 2^N e^{-2\gamma\lambda}. \tag{5.1}$$

By Proposition 2.2,

$$D_t - D_s = -B([s,t]) + \mathbb{E}A([0,s];[s,t]) - A([0,s];[s,t])$$
$$\leq -B([s,t]) + c_1(t-s)L(t-s).$$

Since $\lambda/(400j^2) > 2c_1(t-s)L(t-s)$, then for $D_t - D_s$ to be larger than $\lambda/(200j^2)$, we must have $-B([t-s]) > \lambda/(400j^2)$. Since $-B([s,t])/(t-s)$ is equal in law to $D_1$, then

$$\mathbb{P}(-B([t-s]) > \lambda/(400j^2)) \leq c_5 \exp(-\gamma\lambda 2^j L(2^j)/(800j^2)).$$

Since for each $j$ there are $2^j$ pairs $(s,t)$ to consider, the probability of (b) is bounded by

$$\sum_{j=N}^{\infty} c_5 2^j \exp(-\gamma\lambda 2^j L(2^j)/(800j^2)).$$

This is summable and can be bounded by

$$c_6 \exp(-\gamma\lambda 2^N/(1600N^3))$$

for some $c_6$. By our choice of $N$, this is less than

$$c_6 e^{-2\gamma\lambda}. \tag{5.2}$$

Combining (5.1) and (5.2), we have

$$\mathbb{P}\left(\sup_{s\leq 1} D_s > \lambda\right) \leq c_7 e^{-2\gamma\lambda}$$

if $\lambda > c_3$. Our result follows immediately from this. $\square$

THEOREM 5.2. *With probability 1,*

$$\limsup_{t\to\infty} \frac{D_t}{t\log\log\log t} \leq \frac{1}{2\pi}.$$

Theorems 5.2 and 5.5 together are just a reformulation of Theorem 1.5.

PROOF OF THEOREM 5.2. Let

$$K = [\log\log t], \quad R = t/K \quad \text{and} \quad I_j = [(j-1)R, jR].$$

Let

$$E_j = \sup_{(j-1)R \leq t \leq jR} (-B([(j-1)R,t])).$$

By Proposition 2.2, if $s < t$ and $(\ell - 1)R \le s < \ell R$, then

$$D_s \le \sum_{j < \ell} (-B(I_j)) + (-B([(\ell-1)R, s]))$$

$$+ \sum_{i < j < \ell} \mathbb{E} A(I_i; I_j) + \sum_{j < \ell} \mathbb{E} A([(\ell-1)R, s]; I_j)$$

$$\le \sum_{j=1}^{K} E_j + \sum_{i < j \le K} \mathbb{E} A(I_i; I_j)$$

$$= \sum_{j=1}^{K} E_j + \sum_{j=1}^{K} \mathbb{E} A([0, (j-1)R]; I_j).$$

By Proposition 2.3 and Remark 2.1, the last term on the last line is bounded by

$$\sum_{j=1}^{K} \frac{1}{2\pi} [jR \log(jR) - (j-1)R \log((j-1)R) - R \log R],$$

which is easily seen to equal

$$\frac{1}{2\pi} t \log \log \log t.$$

Then for $\varepsilon > 0$ and $t$ large enough, we have

$$\mathbb{P}\left(\sup_{s \le t} D_s > (1 + 2\varepsilon) \frac{1}{2\pi} t \log \log \log t\right)$$

$$\le \mathbb{P}\left(\sum_{j=1}^{K} E_j > \varepsilon t \log \log \log t\right)$$

$$= \mathbb{P}\left(\sum_{j=1}^{K} \frac{E_j}{R} > \frac{\varepsilon t \log K}{R}\right) = \mathbb{P}\left(\sum_{j=1}^{K} \frac{E_j}{R} > \varepsilon K \log K\right)$$

$$\le c_6 e^{-\varepsilon K \log K} \mathbb{E} \exp\left(\sum_{j=1}^{K} \frac{E_j}{R}\right) = c_6 e^{-\varepsilon K \log K} \left(\mathbb{E} \exp\left(\frac{E_1}{R}\right)\right)^K,$$

using the independence of the $E_j$. Since $E_1/R$ is equal in law to $\sup_{s \le 1} D_s$, then by Proposition 5.1, the above is bounded by

$$c_6 e^{-\varepsilon K \log K} (c_7)^K.$$

If we take $t$ large enough, we have the bound

$$c_8 e^{-2K}.$$

We apply this with $t_n = q^n$ with $q > 1$ close to 1 so that $(1+3\varepsilon)/q > 1+2\varepsilon$. Since $\exp(-2\log\log t_n) = O(n^{-2})$, we have

$$\mathbb{P}\left(\sup_{s \leq t_n} D_s > (1+3\varepsilon)\frac{1}{2\pi}t_n \log\log\log t_n\right) \leq \frac{c_9}{n^2}$$

for $n$ large. If $D_s > (1+4\varepsilon)\frac{1}{2\pi}s\log\log\log s$ for some $s \in [t_{n-1}, t_n]$, then it follows that $\sup_{s \leq t_n} D_s > \frac{1+3\varepsilon}{q}\frac{1}{2\pi}t_n \log\log\log t_n$. By the Borel–Cantelli lemma, it follows that

$$\limsup_{t \to \infty} \frac{D_t}{t \log\log\log t} \leq \frac{1}{2\pi}(1+4\varepsilon) \quad \text{a.s.}$$

Since $\varepsilon$ is arbitrary, our result follows. $\square$

We now turn to the lower bound.

LEMMA 5.3. *The quantity*

$$\mathbb{P}\left(\sup_{\substack{x \in B(0,3) \\ r < 2}} \frac{1}{r}\int_0^1 \mathbb{1}_{B(x,r)}(X_s)\,ds > \lambda\right)$$

*tends to 0 as* $\lambda \to \infty$.

PROOF. Let

$$U_t(x,r) = \int_0^t \mathbb{1}_{B(x,r)}(X_s)\,ds.$$

By symmetry, the expectation of $\mathbb{E}^y U_1(x,r)$ is largest when $y = x$. We have

$$\mathbb{E}^x \int_0^1 \mathbb{1}_{B(x,r)}(X_s)\,ds = \int_0^1 \int_{B(x,r)} \frac{1}{2\pi s} \exp\left(\frac{-|z-x|^2}{2s}\right) dz\,ds,$$

and a straightforward calculation shows this is bounded by $c_1 r^2(1+\log^+(1/r))$, where $c_1$ can be chosen to be independent of $x$ and $r$. Then by the Markov property,

$$\mathbb{E}[U_1(x,r) - U_t(x,r)|\mathcal{F}_t] = \mathbb{E}^{X_t} U_{1-t}(x,r) \leq c_1 r^2(1+\log^+(1/r)).$$

By [1], Theorem I.6.11, since $U_t(x,r)$ has continuous paths and is nondecreasing, there exists $c_2$ such that

(5.3) $$\mathbb{E}\exp(c_2 U_1(x,r)/r^2(1+\log^+(1/r))) \leq 2.$$

Set $r_k = 2^{-k}$ and let $\mathcal{A}_k$ be the set of points in $B(0,4)$ such that each coordinate is an integer multiple of $2^{-k}$. The cardinality of $\mathcal{A}_k$ is less than $c_3 2^{2k}$. By Chebyshev's inequality,

$$\mathbb{P}\left(\sup_{x_k \in \mathcal{A}_k} \frac{1}{r_k} U_1(x, r_k) > \frac{\lambda}{4}\right) \leq c_4 2^{2k} \exp\left(\frac{-c_5 \lambda}{r_k(1+\log^+(1/r_k))}\right).$$

This is summable in $k$, so

$$\mathbb{P}\left(\sup_{k\geq -1}\sup_{x_k\in\mathcal{A}_k}\frac{1}{r_k}U_1(x,r_k)>\frac{\lambda}{4}\right)$$

tends to 0 as $\lambda\to 0$. If $x\in B(0,3)$ and $r<2$, then $B(x,r)\subset B(x_k,r_k)$ for some $x_k\in\mathcal{A}_k$ and some $k$ such that $r_k/4\leq r\leq r_k$. Our result now follows. $\square$

LEMMA 5.4. *Suppose $\mu$ is a measure supported in $B(0,2)$ such that for all $r\leq 2$ and all $x\in\mathbb{R}^2$, we have $\mu(B(x,r))\leq c_1 r$. There exists $c_2$ such that for all $x\in\mathbb{R}^2$,*

$$\int_0^1\int p(s,x,y)\mu(dy)\,ds\leq c_1 c_2,$$

*where $p(s,x,y)=1/(2\pi s)\exp(-|x-y|^2/(2s))$ is the transition density of two-dimensional Brownian motion.*

PROOF (cf. [2]). Substituting $t=|x-y|^2/(2s)$ shows that

$$\int_0^1 p(s,x,y)\,ds\leq c_3(1+\log^+(1/|x-y|)).$$

We then use the Fubini theorem to write

$$\int_0^1\int p(s,x,y)\mu(dy)\,ds$$

$$\leq c_4\sum_{k=-1}^\infty\int_{B(x,2^{-k})\setminus B(x,2^{-k-1})}(1+\log^+(1/|x-y|))\mu(dy)$$

$$\leq c_5\sum_{k=-1}^\infty(2+k)\mu(B(x,2^{-k}))\leq c_6 c_1\sum_{k=-1}^\infty(2+k)2^{-k}\leq c_7 c_1,$$

as required. $\square$

THEOREM 5.5. *We have*

$$\limsup_{t\to\infty}\frac{-\beta_t}{t\log\log\log t}\geq\frac{1}{2\pi}\qquad a.s.$$

PROOF. Let $K=[b\log\log t]$ and $R=t/K$, where $b$ is to be chosen later. Let $I_j=[(j-1)R,jR]$. Let $\mathcal{G}_j=\sigma(X_s:s\leq jR)$.

By (2.5) we have

$$-\beta_t=\sum_{j=1}^K -B(I_j)-\sum_{j=1}^K A(I_j;[0,(j-1)R])+\sum_{j=1}^K \mathbb{E}A(I_j;[0,(j-1)R])$$

$$=J_1+J_2+J_3.$$

Recall that $A(I_j; [0, (j-1)R])$ is equal in law to $\alpha(R, (j-1)R)$. By Proposition 2.3 and Remark 2.1,

$$J_3 = \sum_{j=1}^{K} \frac{1}{2\pi}[jR\log(jR) - R\log R - (j-1)R\log((j-1)R)]$$

and it is straightforward to see that this is equal to $\frac{1}{2\pi}t\log K$.

Define the sets

$D_{j1} = \{X_{jR} \in B(j\sqrt{R}, \sqrt{R}/16)\}$,

$D_{j2} = \{X(I_j) \subset [(j-1)\sqrt{R} - (\sqrt{R}/8), j\sqrt{R} + (\sqrt{R}/8)] \times [-\sqrt{R}/8, \sqrt{R}/8]\}$,

$D_{j3} = \{B(I_j) \leq \kappa_1 R\}$,

$D_{j4} = \Big\{ \int_{(j-1)R}^{jR} \mathbb{1}_{B(x,r\sqrt{R})}(X_s)\, ds \leq \kappa_2 rR$

$$\text{for all } x \in B(j\sqrt{R}, 3\sqrt{R}), 0 < r < 2\sqrt{R}\Big\},$$

$D_{j5} = \{A(I_{j-1}; I_j) \leq \kappa_3 R\}$,

where $\kappa_1, \kappa_2$ and $\kappa_3$ are constants to be chosen later and that do not depend on $j, b, t$ and $R$. Let

$$C_j = D_{j1} \cap D_{j2} \cap D_{j3} \cap D_{j4} \cap D_{j5}$$

and

$$E = \bigcap_{j=1}^{K} C_j.$$

We want to show

(5.4) $$\mathbb{P}(C_j | \mathcal{G}_{j-1}) \geq c_1$$

on the set $C_1 \cap \cdots \cap C_{j-1}$, where $c_1 > 0$ does not depend on $j, b, t$ and $R$. Once we have (5.4), then

$$\mathbb{P}\Big(\bigcap_{i=1}^{m} C_i\Big) = \mathbb{E}\Big(\mathbb{P}(C_m | \mathcal{G}_{m-1}); \bigcap_{i=1}^{m-1} C_i\Big) \geq c_1 \mathbb{P}\Big(\bigcap_{i=1}^{m-1} C_i\Big)$$

and, by induction,

(5.5) $$\mathbb{P}\Big(\bigcap_{i=1}^{K} C_i\Big) \geq c_1^K = c_1^{b\log\log t} = \exp(b\log\log t \log c_1).$$

On the set $E$ we have

$$-J_1 = \sum_{j=1}^{K} B(I_j) \leq \kappa_1 KR = \kappa_1 t.$$

SELF-INTERSECTION LOCAL TIME 25Since for each $j$ we are on the set $D_{j2}$, then on the event $E$ we have $X(I_i)$ disjoint from $X(I_j)$ if $|i-j| > 1$. Therefore,

$$-J_2 = \sum_{j=1}^{K} A(I_{j-1}; I_j) \leq \kappa_3 KR = \kappa_3 t.$$

So on $E$ we have altogether that

$$J_1 + J_2 + J_3 \geq \frac{1}{2\pi} t \log\log\log t - (\kappa_1 + \kappa_3)t.$$

We now proceed to show (5.4). By the support theorem for planar Brownian motion and scaling (see [1], Theorem I.6.6),

$$\mathbb{P}(D_{j1} \cap D_{j2} | \mathcal{G}_{j-1}) > c_2.$$

By scaling,

$$\mathbb{E}B(I_j) = R\mathbb{E}B([0,1]) \leq c_3 R.$$

So if $\kappa_1$ is chosen large enough,

$$\mathbb{P}(B(I_j) > \kappa_1 R) < c_2/6.$$

Now let us look at $D_{j4}$. By scaling and Lemma 5.3 it follows that $\mathbb{P}(D_{j4}) \leq c_2/6$ if we choose $\kappa_2$ large enough.

Next we look at $D_{j5}$. We have an estimate on $\mathbb{E}A(I_{j-1}; I_j)$, but what we actually need is an estimate on $\mathbb{E}[A(I_{j-1}; I_j)|\mathcal{G}_{j-1}]$. To show

$$\mathbb{P}(A(I_{j-1}; I_j) > \kappa_3 R | \mathcal{G}_{j-1}) \leq c_2/6$$

if $\kappa_3$ is large enough, it is enough to show

(5.6) $$\mathbb{E}[A(I_{j-1}; I_j)|\mathcal{G}_{j-1}] \leq c_4 R$$

on the set $\bigcap_{i=1}^{j-1} C_i$. By [3], we can let $\mu$ be the measure on $\mathbb{R}^2$ defined by

$$\mu(F) = \int_0^1 \mathbb{1}_F(X_s)\, ds$$

and consider $A([0,1]; [1,2])$ as an additive functional of Brownian motion that corresponds to the measure $\mu$. So to show (5.6), by scaling and translation invariance it is enough to show

(5.7) $$\mathbb{E}[A([0,1]; [1,2])|\mathcal{F}_1] \leq c_5$$

on the set where $\mu(B(x,r)) \leq \kappa_1 r$ for all $x \in B(0,3)$ and all $r \in (0,2)$. The conditional expectation (5.7) is bounded by

$$\sup_y \int_0^1 \int p(s, y, z) \mu(dz)\, ds,$$

where $p(s, y, z)$ is the transition density for planar Brownian motion. Using Lemma 5.4 we have (5.7).

Setting $c_6 = c_2/2$ gives the desired lower bound (5.4).

Let $t_n = \exp(n^\gamma)$ for some $\gamma > 1$ and let $\varepsilon > 0$. Provided we take $b$ (in the definition of $K$) small enough, the Borel–Cantelli lemma tells us that

$$(5.8) \quad -B([t_{n-1}, t_n]) \geq \left(\frac{1}{2\pi} - \varepsilon\right)(t_n - t_{n-1}) \log \log \log(t_n - t_{n-1})$$

infinitely often with probability 1. By (2.5) we have

$$(5.9) \quad \begin{aligned} -B([0, t_n]) =& -B([0, t_{n-1}]) - B([t_{n-1}, t_n]) - A([0, t_{n-1}]; [t_{n-1}, t_n]) \\ &+ \mathbb{E}A([0, t_{n-1}]; [t_{n-1}, t_n]). \end{aligned}$$

By the upper bound for $-\beta_t$, from Theorem 4.3, we know that

$$(5.10) \quad B([0, t_{n-1}]) = O(t_{n-1} \log \log \log t_{n-1}) = o(t_n \log \log \log t_n) \quad \text{a.s.}$$

By scaling,

$$(5.11) \quad \widetilde{A}_n = A([0, t_{n-1}]; [t_{n-1}, t_n])$$

is equal in law to $\Delta_n \alpha(1, t_{n-1}/\Delta_n)$, where $\Delta_n = t_n - t_{n-1}$. By Lemma 4.1,

$$\mathbb{P}(\widetilde{A}_n > t_n) \leq \exp\left(-c_8 \frac{t_n}{\Delta_n} \sqrt{\frac{\Delta_n}{t_{n-1}}}\right),$$

which is summable. Using the Borel–Cantelli lemma, we have

$$(5.12) \quad \widetilde{A}_n = o(t_n \log \log \log t_n) \quad \text{a.s.}$$

Substituting this, (5.8), (5.10) and (5.12) in (5.9) proves the theorem. $\square$

REMARK 5.6. Theorem 1.5 follows immediately from Theorems 5.4 and 5.5.

The following corollary completes the proof of Theorem 1.3.

COROLLARY 5.7. *Let $L$ be as in the statement of Theorem* 3.2. *Then $L < \infty$.*

PROOF. In the proof of Theorem 5.5 we showed that the event $E$ had probability at least $\exp(b \log \log t \log c_1)$ and that on the event $E$ we had $-\beta_t \geq \frac{1}{2\pi} t \log \log \log t - c_2 t$ provided $t$ was large enough. Choose $t$ so that $\frac{1}{2\pi} \log \log \log t - c_2 = \log s$. Using scaling, we then have

$$\mathbb{P}(-\beta_1 > \log s) \geq \exp(-c_3 \log \log t) = \exp(-c_4 s^{2\pi}).$$

Now take logarithms of both sides and divide by $s^{2\pi}$. $\square$

REMARK 5.8. Theorem 1.3 follows immediately from Theorem 3.2 and Corollary 5.7.

**Acknowledgments.** We thank A. Dorogovtsev, J.-F. Le Gall, W. Li, J. Rosen and Z. Shi for their interest in this work and for helpful discussions. We also thank D. Bakry, E. Carlen, J. Denzler and C. Villani for helpful information about the best constants for Gagliardo–Nirenberg inequalities. Finally, we thank the referee and an Associate Editor for their exceptionally helpful reports.

DEPARTMENT OF MATHEMATICS
UNIVERSITY OF CONNECTICUT
STORRS, CONNECTICUT 06269-3009
USA
E-MAIL: bass@math.uconn.edu

DEPARTMENT OF MATHEMATICS
UNIVERSITY OF TENNESSEE
KNOXVILLE, TENNESSEE 37996-1300
USA
E-MAIL: xchen@math.utk.edu